\documentclass[twocolumn]{autart}    
\usepackage{tikz}
\usepackage{tikz-3dplot}
\usepackage{pgfplots}
\usepgfplotslibrary{groupplots} 
\usetikzlibrary{pgfplots.groupplots} 

\usepackage{url}
\usepackage{amsmath,amssymb,amsfonts,stmaryrd}
\usepackage{graphicx}          
\usepackage{color}
\newtheorem{definition}{\bf Definition}[section]

\newtheorem{remark}{\bf Remark}[section]
\newtheorem{proposition}{\bf Proposition}[section]
\newtheorem{theorem}{\bf Theorem}

\newcommand{\ra}[1]{{\color{black}#1}} 
\newcommand{\rc}[1]{{\color{black}#1}} 

\begin{document}

\begin{frontmatter}

    \title{Large Deviations in Safety-Critical Systems with Probabilistic Initial Conditions} 

\thanks[footnoteinfo]{This paper was not presented at any IFAC
meeting. Corresponding author Aitor R. Gomez.}


\author[Paestum]{Aitor R. Gomez}\ead{arg@es.aau.dk},    
\author[Rome]{Manuela L. Bujorianu}\ead{l.bujorianu@ucl.ac.uk}, 
\author[Paestum]{Rafal Wisniewski}\ead{raf@es.aau.dk}  

\address[Paestum]{Section of Automation \& Control, Aalborg University, DK}  
\address[Rome]{Department of Computer Science, University College London, UK}             

\begin{keyword}                           
safety; rare events; optimal control problem; random initial conditions; large deviations; approximation.               
\end{keyword}                             

\begin{abstract}                          
We often rely on probabilistic measures---e.g. event probability or expected time---to characterize systems' safety. However, determining these quantities for extremely low-probability events is generally challenging, as standard safety methods usually struggle due to conservativeness, high-dimension scalability, tractability or numerical limitations. We address these issues by leveraging rigorous approximations grounded in the principles of Large Deviations theory. By assuming deterministic initial conditions, Large Deviations identifies a single dominant path in the low-noise limit as the most significant contributor to the rare-event probability: the \emph{instanton}. We extend this result to incorporate stochastic uncertainty in the initial states, which is a common assumption in many applications. To that end, we determine an expression for the probability density of the initial states, conditioned on the unsafe rare event being observed. This expression gives access to the most probable initial conditions, as well as the most probable hitting time and path deviations, leading to the realization of the unsafe event. We demonstrate it's effectiveness by solving a high-dimensional and non-linear problem: a space collision.

\end{abstract}

\end{frontmatter}

\section{Introduction}\label{sec:introduction}

A complete mathematical characterization of the physical interactions between a system and its environment is generally intractable. Inevitably, uncertainties---although often small in magnitude---persist, and can gradually influence the long-term behavior of the system, leading to deviations from a reference domain.
Examples of this are abundant in areas such as spacecraft collision avoidance, air traffic management, missile guidance, telescope and radio pointing, and others \cite{Meerkov:88}.

The notion of safety in control engineering is usually associated with critical constraints that can not be violated for the sake of the systems sustainability. Viability kernels \cite{Aubin:09} and minimal (or maximal) reachability analysis \cite{Ian:07} are solid frameworks that have been successfully applied to prove and preserve safety. Both provide in-depth insights upon non-deterministic \cite{Ian:12} and stochastic \cite{Bujorianu:12} systems' trajectories by identifying sets of initial (or final) states for which the system is, in some formal sense, safe. These two concepts are key to guarantee safety, yet their formulations rely on set-valued time propagation of the system dynamics through partial differential equations (PDEs). Exact solutions to these problems are generally very hard---or even impossible---to determine, and numerical approximations are prohibitive when exceeding more than four state dimensions due to Bellman's \emph{curse of dimensionality}.\footnote{This is an open area of research and such considerations are made at the current time of writing this report.} Such limitations restrict these formulations mostly to systems with low state space dimensions and/or linear dynamics.


Barrier certificates \cite{Prajna:04} emerged from the desire of verifying systems safety while avoiding the difficult computation of reachable sets or density propagation.
Certificates are based on showing the existence of a function with zero-level set separating the trajectories of the system issued from a given initial set, from the unsafe region of the state space. When the uncertainties and disturbances are of stochastic nature and have no hard bounds, it is more reasonable to define safety from a probabilistic standpoint. A common approach within this framework is using supermartingale theory to construct a stochastic counterpart of a barrier that yields an upper bound on the probability of reaching the unsafe region. However, stochastic barriers, and the closely related concept of \emph{p}-safety in reach-avoid problems \cite{Rafal:21,Manuela:20}, involve a Fokker-Planck PDE that makes the problem less tractable, thereby requiring several relaxations to enable its computation. Stochastic barrier techniques have been attempted in a previous work \cite{Aitor:22}, which demonstrated various accuracy and numerical limitations due to such relaxations and the extremely low probabilities describing the event.

On the other hand, estimating the probability of rare events using crude Monte Carlo or particle methods is generally simpler but impractical, as it requires an infeasible number of samples to achieve significant estimates. Furthermore, including uncertainty in the initial conditions in addition to process uncertainty, as contemplated in this paper, exacerbates this challenge even further.

In order to overcome these limitations---tractability, high-dimensional scalability, numerical instabilities, and accuracy---we leverage results from Large Deviations theory (LDT). LDT provides an asymptotic approximation of the system's limiting behavior: as the process noise decreases, the probability density becomes increasingly concentrated around the unperturbed (deterministic) path. Consequently, the probability contribution of all system's trajectories deviating away from the deterministic path to reach an unsafe region decays exponentially, leaving only the contribution of a single dominant path: the \emph{instanton}, or maximum likelihood path. This dominant path is found as the solution to a variational problem (VP) with fixed initial conditions and final time. Usually, however, the initial conditions and the (final) time at which the unsafe event takes place are random. Depending on the system's dynamics and characterization of the unsafe region, random final times might result in the instanton not being a unique solution; this fact is already accounted for in LDT. On the contrary, LDT lacks formulations accounting for non-deterministic initial conditions, e.g. when they are described stochastically.

Approximating stochastic processes with LDT is not novel; it has been widely employed in various optimal control problems in the past \cite{Meerkov:88,Runolfsson:92,Dupuis:87,Whittle:90,Matthias:19}. The well-known work of Freidlin and Wentzell \cite{Freidlin:84} established the theory of Large Deviations for stochastic differential equations, forming the basis for the main results of this paper. Kushner later extended the results for degenerate diffusion processes \cite{Kushner:84,Agnessa:10}, which aligns more accurately with the description of mechanical systems. While the theory can be abstract, some works \cite{Tobias:13,Matthias:19,Touchette:09} offer a more accessible introduction tailored for engineers.


This collection of results on LDT provide a convenient, rigorous and precise mathematical framework that can be used for estimating probabilistic quantities of rare events more accurately when the initial condition is deterministic. Note that, by relying on the maximum principle to solve the VPs, the issue of dealing with high dimensional state spaces can be alleviated, as opposed to using the complex viscosity solution techniques employed in reachability and Fokker-Planck PDEs. The method is numerically more stable than stochastic barriers, as it characterizes the probability as an integral over a path, instead of zero level-sets of functions; it also gives access to various quantities of interest, e.g. expected time of occurrence or most likely deviations; and it can be used in combination with other techniques, e.g. importance sampling, to perform simulations of rare events. Additionally, there are readily available numerical tools to solve VPs, allowing its seamless implementation in a wide range of applications.

\subsection{Organization and Contributions}

The necessary background on LDT and the weak \emph{p}-safety problem are introduced in section \ref{sec:preliminaries}. Section \ref{sec:uncertainty} constitutes the main contributions: we approximate the posterior probability density of the initial states conditioned on the observation of a rare event, and provide a computational expression for the solution of the weak \emph{p}-safety problem in weakly-perturbed systems. 
Section \ref{sec:example} is devoted to a real, non-linear and high-dimensional example involving spacecraft collisions. Necessary conditions for optimality are derived using the maximum principle, which are used to cross-check the obtained numerical solution. Finally, in section \ref{sec:conclusions} we give the final remarks and conclude the paper. The general contributions of the paper are the following: 
\begin{enumerate}
    \item[i)] Extended methodology for safety analysis of dynamical systems to include unsafe rare events.
    \item[ii)] Theorem \ref{thm:main}, enabling an expression for the initial state distribution describing how probable is to depart from a point given an unsafe region is visited.
\end{enumerate}
Various safety-related questions can be addressed using Theorem \ref{thm:main}, either directly through its solution or as an intermediate step toward more advanced analyses. For instance, in the context of an in-orbit satellite collision, it allows us to estimate the most probable collision configurations---trajectory and time of collision---which can improve greatly the subsequent studies and tracking of the resulting debris. We shall emphasize that the resulting expression in Theorem \ref{thm:main} is not analytical; thus, the calculation of the total probability of collision (i.e. the solution to weak \emph{p}-safety) becomes non-trivial. This is a direct consequence of the inclusion of uncertainty in the initial state. Nonetheless, Theorem \ref{thm:main} plays a crucial role in enabling methods, such as importance sampling or Bayesian optimization, which offer a computationally tractable approach to accurately approximate the associated integral. Additionally, Theorem \ref{thm:main} can be extended as a min-max optimization problem to determine control inputs that minimize the likelihood of a collision.

\subsection{Notation}
We will consider systems in $\mathbb{R}^{n}$, with $\mathcal{C}^0(I,\mathbb{R}^{n})$ denoting the space of continuous functions that maps the time interval $I$ into the state space. For mechanical systems we consider a space of position and momenta $\mathbb{R}^{m}\!\times\!\mathbb{R}^{m}$, and use $\mathcal{C}^0(I,\mathbb{R}^{m}\!\times\!\mathbb{R}^{m})$. Occasionally, we will refer to it as $\mathcal{C}_I$ to simplify notation. We will assume a probability space $(\Omega,\mathcal{F},\textbf{P})$, identifying $\Omega$ with the space of trajectories $\mathcal{C}_I$. We write $X_t(\omega)$ and $X(\omega,t)$ interchangeably, denoting a point in a sample path of the process. We also use $\phi$ to denote a sample path $X^\epsilon(\omega)$ in the limit $\epsilon\rightarrow0$, and $\phi_{y}$ to specify that $\phi(0)=y$. $\textbf{P}_{\mu}$ refers to a probability with respect to a measure $\mu$, and $\textbf{P}_{y}$ a probability conditioned on the initial point $y$. We often use the notation $x=(y,z)\equiv[y^\top,z^\top]^\top$ to gather column vectors $y,z\in\mathbb{R}^n$ into a column vector $x\in\mathbb{R}^{2n}$. For a vector $v$, the 2-norm is denoted as $||v||$, and the $\infty$-norm as $||v||_\infty$. The notation $\nabla_x h$ is used to refer to the gradient of a function $h$ with respect to $x$. For a set $D$ subset of a topological space, its boundary is $\partial D$ and its closure $\text{cl}(D)$.

\section{Problem formulation and Background}\label{sec:preliminaries}
\subsection{System Description}
We will assume that unmodeled external perturbations are Gaussian and small in nature, and describe the system by stochastic differential equations (SDE) in $\mathbb{R}^n$, emphasizing a small noise parameter $0<\epsilon\ll1$
\begin{equation}\label{eq:dSDE_X}
    dX_t^\epsilon= b(X_t^\epsilon)dt+\sqrt{\epsilon}\sigma dW_t,\quad \rc{X^\epsilon(0)\sim \mathcal{N}(x_0,\Sigma)},
\end{equation}
where $b:\mathbb{R}^n\rightarrow\mathbb{R}^n$ is Lipschitz continuous and denotes modelled forces driving the system, $\sigma\in\mathbb{R}^{n\times n}$ is the diffusion matrix defining the noise covariance $a=\sigma\sigma^\top$, and $(W_t)$ is an $n$-dimensional Wiener process. We regard $(X_t^\epsilon)$ as the process resulting from weakly perturbing the deterministic system (i.e. for $\epsilon=0$)
\begin{equation}\label{eq:ode_deterministic}
    \dot{x}=b(x),\quad x(0)=x_0.\\[-5mm]
\end{equation}

This work is motivated by the search for a suitable safety framework that handles low-probability events particularly in mechanical systems, for which we may replace \eqref{eq:dSDE_X} with the following second-order SDE
\begin{equation}\label{eq:dSDE}
    \begin{aligned}
        dR_t =& V_t dt,\\
        dV_t =& b(R_t)dt+ \sqrt{\epsilon}\sigma dW_t,
    \end{aligned}
\end{equation}
where $R(t)$ and $V(t)$ are random position and momentum vectors respectively, in $\mathbb{R}^m$. Similarly, we regard the previous SDE as the perturbed version of $\ddot{r}=b(r)$.

We will refer when necessary to the nuances arising when replacing \eqref{eq:dSDE_X} with \eqref{eq:dSDE}.

 \subsection{Safety Problem}
Our aim is two-fold: 1) to find an expression that approximates the unsafe distribution of initial states, which describes how probable is that $(X_t^\epsilon)$ departs from \rc{a point sampled by a distribution $\mu$ and subsequently visits an unsafe compact set $D\subset\mathbb{R}^n$}, and 2) to find the most probable configurations (i.e. initial states, hitting times and path deviations) in which this event is observed. Events of this type will be characterized by very low probabilities, consequence of the smallness of $\epsilon$, earning thus the nomination of Unsafe Rare Event (URE). Many practical interpretations can be given to $D$. For instance, a region of space occupied by a piece of debris on a spacecrafts orbit, undesired targets for a missile or other vehicles near aircraft pathways.

The problem of a random system visiting a set $D$ can be defined using random stopping times, which motivates the following definition.
\begin{definition}[Hitting time]\label{def:hitting}
     Consider the system \eqref{eq:dSDE_X} hitting the boundary of the set $D$ at most once at time
     $$\tau_D^\epsilon(\omega)=\inf\{ t>0 : X^\epsilon(\omega, t)\in \partial D \}.$$\\[-12mm]
 \end{definition}
 This is the \emph{first} hitting time of the boundary which is equal to the first hitting time of $D$, since the trajectories of the process are continuous.
 \begin{remark}
     We may have an URE associated only to some components of $X^\epsilon(t)$. This is common in mechanical systems, where the event could be represented as $\{R(\omega,t)\in\partial \tilde{D}\}$;\footnote{In mechanical systems, defining $\tilde{D}$ as a compact set only on $\mathbb{R}^m$ is not a restriction since we can define an unsafe set for $V(t)$ to be large enough.} see for instance example in section \ref{sec:example}.
 \end{remark}
\rc{Initially, we strive to compute the probability that $(X_t^\epsilon)$ reaches $\partial D$ in a finite amount of time $T$, assuming $X_0^\epsilon$ is random. For an initial probability measure $\mu$, we write this in the following familiar form
\begin{equation}\label{eq:psafety}
    \textbf{P}_\mu[\tau_D^\epsilon\leq T]=\int\textbf{P}_{y}[\tau_D^\epsilon\leq T]\mu(dy),
\end{equation}
which coincides with the definition of weak $p$-safety \cite{Rafal:21}. Note that the hitting event $E_y\equiv\{\omega\!:\!\tau_D^\epsilon(\omega)\!\leq\!T|X_0^\epsilon=y\}$ corresponds to all the (unsafe) paths departing from $y$, which we gather in a set $H_D(y)\!\subset\mathcal{C}_I$ defined as follows.
\begin{definition}[Set of unsafe paths]
    Given a point $y$, a time $T$ and a (compact) region $D$, the unsafe paths are defined as all continuous functions issued from $y$ and reaching $\partial D$ at some time $t\in I\equiv[0,T]$. That is,
        \begin{equation*}
            H_D(y)\equiv\!\{\phi\in\mathcal{C}_I : \phi(0)=y, ~\exists t\in[0,T]~ \text{ s.t. }~\phi(t)\in\partial D\}.
        \end{equation*}
\end{definition}
The large deviations theory introduced by Freidlin and Wentzell \cite{Freidlin:84} provides a formal way of treating the probability of the hitting event $E_y$ in the limiting case $\epsilon\rightarrow0$. This implies that accurate approximations can be leveraged employing the same techniques, provided the process noise is small enough. We review the necessary tools to establish these relationships in the next section.
}

\subsection{Large Deviations Approximation}

The asymptotic solutions developed for sufficiently small $\epsilon$ are well known in the mathematical and physics community (see e.g. \cite{Freidlin:84,Touchette:09,Touchette:13,Dembo:10} and the references within), which provide formal approximations of path probabilities under certain regularity conditions. In order to introduce these methods, let us briefly analyse the limiting behavior of the weakly perturbed system \eqref{eq:dSDE_X} as a function of the noise parameter. For now, consider that $X^\epsilon(0)=x_0$ with probability 1. Then, note that as we decrease $\epsilon$, the process $(X_t^\epsilon)$ is expected to remain closer to the deterministic trajectory $x$ (the solution of \eqref{eq:ode_deterministic}). In the limit $\epsilon\rightarrow0$, we can formally write that, for $\delta>0$
\begin{equation}\label{eq:decay}
\rc{
    \textbf{P}_{x_0}\big\{\omega:d(X^\epsilon(\omega),x)\leq\delta\big\}\rightarrow1,
}
\end{equation}
\rc{which holds uniformly in small $\delta$. The metric $d$ is defined by the sup-norm in the path space  $$d(\psi,\phi)\equiv||\psi-\phi||_\infty=\underset{t\in I}{\sup}~||\psi_t-\phi_t||.$$}
The expression in \eqref{eq:decay} reveals that the probability of fluctuating away from the deterministic solution $x$ decays to zero with decreasing $\epsilon$. In fact, this decay in probability allows us to \ra{precisely characterize the probability of $(X_t^\epsilon)$ being close to a given path $\phi\in\mathcal{C}_I$ issued from $x_0$, in the limit $\epsilon\!\rightarrow\!0$. This is a key result known as Large Deviation Principle (LDP) for diffusions.
\begin{equation}
    \lim_{\epsilon\rightarrow 0} \epsilon\ln \textbf{P}_{x_0}\big[d(X^\epsilon,\phi)\leq\delta\big]=-S_T(\phi),
    \label{eq:LDP_1}
\end{equation}}where $S_T:\mathcal{C}_I\rightarrow\mathbb{R}_{\geq 0}$ is referred to as the (normalized) action functional, rate function or path entropy (depending on the field of study) on the interval $I\equiv[0,T]$. We will simply call it action functional and define it as \cite{Kushner:84,Agnessa:10}
\rc{
\begin{equation}\label{eq:action}
\begin{aligned}
    S_T(\phi) &\equiv\frac{1}{2}\int_0^T||w_t||^2dt,
\end{aligned}
\end{equation}}if $\phi$ is absolutely continuous, and $S_T(\phi)=+\infty$ otherwise. The value of $S_T(\phi)$ represents the total deviation from the vector field $b$ accumulated by the path $\phi$. For the example in section \ref{sec:example}, it may even be plausible to give a physical interpretation to this value. Note that $S_T$ is convex with minimum at $w=0$ (i.e. no perturbations), which yields $\dot{\phi}=b(\phi)$ (or $\dot{\nu}=b(\eta)$ for mechanical systems), making its solution equivalent to that of the deterministic system $\dot{x}=b(x)$ (or $\ddot{r}=b(r)$). This implies that the deterministic trajectory is the most probable path the system will follow, i.e. it exhibits highest probability, in the sense of the LDP.

In the case of extremely low process noise, we may regard the sample paths of the process $(X_t^\epsilon)$, as absolutely continuous functions solving the differential equation
\begin{equation}\label{eq:phi_dynamics}
    \dot{\phi}=b(\phi)+\sigma w,\quad \phi(0)=x_0,
\end{equation}
where $w:I\rightarrow\mathbb{R}^n$ is a measurable function representing deviations from the vector field $b$, which steer the system \eqref{eq:phi_dynamics} mimicking the perturbations of the real system \eqref{eq:dSDE_X}.

\begin{remark}
    In the mechanical system description, we shall specify that $\phi(t)=(\eta(t),\nu(t))$ where $\eta$ and $\nu$ are continuous functions representing sample paths of $(R_t)$ and $(V_t)$ respectively. The governing equations are then
    \begin{equation}\label{eq:phi_dynamics_mech}
        \begin{aligned}
            \dot{\eta}&=\nu,\\
            \dot{\nu}&=b(\eta)+\sigma w.
        \end{aligned}
    \end{equation}
\end{remark}
\vspace{-0.5cm}
The involvement of the action functional $S_T$ in the probability of $(X_t^\epsilon)$ lying in the neighborhood of one such path incentivizes a refined definition of the unsafe set.
\rc{
\begin{definition}[Unsafe paths with finite action]
    An unsafe path $\phi\in H_D(x_0)$ is well-behaved if it belongs to the following set
    \begin{equation*}
        \Phi_D(x_0)\equiv\!\{\phi\in H_D(x_0) : ~\text{s.t.}~\eqref{eq:phi_dynamics},~S_T(\phi)<\infty\},
    \end{equation*}
    which is compact in the sup-norm topology \cite{Freidlin:84}.
\end{definition}
}
\ra{
The extreme path $\phi^\ast$ that infimizes $S_T$ on $\Phi_D(x_0)$ serves as a lower bound to the probability of $E_{x_0}$, given that
\begin{equation}\label{eq:lowerb}
    \textbf{P}_{x_0}[\tau_D^\epsilon\leq T]\geq \textbf{P}_{x_0}[d(X^\epsilon,\phi^\ast)\leq\delta],
\end{equation}
which is generally true for small enough $\delta>0$. It is perhaps more remarkable that the reverse inequality can also be proven in the asymptotic case. To build intuition about this fact, consider $\delta$-balls centered around finitely many continuous paths $\{\phi_i\}_{i=1}^n$ covering the closure of $H_D(x_0)$, and include $\phi^\ast$. We may use paths from $\Phi_D(x_0)$ and $\Phi_D^c(x_0)\equiv \text{cl}(H_D)\backslash\Phi_D$. Nevertheless, paths in $\Phi_D^c(x_0)$ have infinite action and thus do not contribute to the probability by definition. Then, for $\phi_i\in\Phi_D(x_0)$
\begin{equation}\label{eq:upperb}
    \textbf{P}_{x_0}[\tau_D^\epsilon\leq T]\leq \sum_{i=1}^n\textbf{P}_{x_0}[d(X^\epsilon,\phi_i)\leq\delta].
\end{equation}
A consequence of \eqref{eq:LDP_1} is that any trajectory of $(X_t^\epsilon)$ contributing to the event $E_{x_0}$ is close to $\phi^\ast$ with overwhelming probability\footnote{Language used in \cite{Freidlin:84} to refer to a very large concentration of the probability distribution.}, while being close to other paths $\phi_i\in\text{cl}(H_D)$ with slightly larger action than $S_T(\phi^\ast)$ contributes negligibly. In the limit, this effectively yields
\begin{equation}\label{eq:asympt}
    \sum_{i=1}^n\textbf{P}_{x_0}[d(X^\epsilon,\phi_i)\leq\delta] \overset{\epsilon}{\asymp} \textbf{P}_{x_0}[d(X^\epsilon,\phi^\ast)\leq\delta].
\end{equation}
Here we used ``$\overset{\epsilon}{\asymp}$" to denote a log-asymptotic equivalence as $\epsilon\!\rightarrow\!0$. Both, \eqref{eq:lowerb} and \eqref{eq:upperb} together with \eqref{eq:asympt}, lead to a sharp result, establishing the exact asymptotic rate of the hitting probability.
\begin{equation}\label{eq:LDP}
    \lim_{\epsilon\rightarrow0}~\epsilon\ln\textbf{P}_{x_0}[\tau_D^\epsilon\leq T]=-\inf_{\phi\in\Phi_D(x_0)}S_T(\phi).
\end{equation}
}This gives rise to the following (well-known) definition.

\begin{definition}[Maximum likelihood pathway]
    The extreme function $\phi^\ast$ infimizing the action functional $S_T$ on $\Phi_D(x_0)$, thus experiencing minimal deviation $w^\ast$, is the maximum likelihood (ML) pathway, or instanton.
\end{definition}
\rc{In practice, we find the instanton by solving a variational problem, the solution of which is given by the pair $(\phi^\ast,t^\ast)$, i.e. the sample path and hitting time, that contributes the most to the probability of $E_{x_0}$. The following is the corresponding variational problem, which is equivalent to \eqref{eq:LDP}.
\begin{equation}\label{eq:quasipotential}
    Q(x_0)\equiv\inf_{t\leq T}\ \underset{\substack{\phi:\phi(0)=x_0,\\ \ \phi(t)\in\partial D}}{\textup{min}} ~S_t(\phi),\quad \text{s.t.}~\eqref{eq:phi_dynamics},
\end{equation}
where $Q$ is a new \emph{action} function depending only on the initial point $x_0$. This function is often referred to as \emph{quasi-potential}\footnote{The usual notation for the quasipotential is $V$. We use $Q$ so as not to confuse it with the random velocity $V(t)$.} (see e.g. \cite{Freidlin:84,Tobias:13}).
}

\section{Randomness in the initial conditions}\label{sec:uncertainty}


In this section, we study the same problem as before but considering the departing state to be random. Specifically, we consider $X^\epsilon(0)\sim \mathcal{N}(x_0,\Sigma)$ and define the initial measure $\mu(dy)$ with density $p_0$ being
\begin{equation}\label{eq:mu}
    p_0(y) = \frac{1}{\sqrt{(2\pi)^{2n}\det \Sigma}}\exp\left(-S_0(y)\right),
\end{equation}
and
\begin{equation}\label{eq:S0}
    S_0(y)\equiv\frac{1}{2}||y-x_0||^2_\Sigma.
\end{equation}
This modification introduces additional degrees of freedom to the problem. Determining the most likely path realizing the URE no longer depends uniquely on the path deviations $w$, but also on the state of departure.

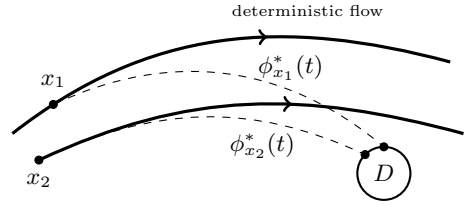
\begin{figure}[ht!]
    \centering
    \usetikzlibrary{decorations.markings}

\begin{tikzpicture}

\newcommand\ang{25}

%

\draw[black, very thick]
    (0,0) to[out=\ang, in=190-\ang] 
    coordinate[pos=0.6] (arrow1A) 
    coordinate[pos=0.61] (arrow1B)
    (160px, 10px);
    
\draw[black, very thick]
    (-10px,10px) to[out=\ang+13, in=190-\ang] 
    coordinate[pos=0.1] (mu) 
    coordinate[pos=0.6] (arrow2A) 
    coordinate[pos=0.61] (arrow2B)
    node[pos=0.7,above=1.5mm] {\tiny deterministic flow}
    (155px, 37px);

\coordinate (D) at (130px,-5px);

\draw[black, thick] (D) circle (10px) node[] {\small $D$};
\draw[black,fill] (mu) circle (1.5pt) node[above=1mm] {\small $y_1$};
\draw[black,fill] (0px,0px) circle (1.5pt) node[below=0.6mm] {\small $y_2$};
\draw[black,fill] (D)++(-7.0711px,7.0711px) circle (1.5pt);
\draw[black,fill] (D)++(0px,10px) circle (1.5pt);
    
\draw[black, dashed]
    (mu) to[out=\ang, in=140] 
    node[pos=0.7, above=0.1mm] {\small $\phi_{y_1}^\ast(t)$} 
    (130px,5px);
    
\draw[black, dashed]
    (0,0) to[out=\ang, in=155] 
    node[pos=0.69, below=1px] {\small $\phi_{y_2}^\ast(t)$} 
    (122.9289px,2.0711px);

\draw[->,very thick] (arrow1A) -- (arrow1B);
\draw[->,very thick] (arrow2A) -- (arrow2B);
\end{tikzpicture}
    \caption{Geometric description of the problem.}
    \label{fig:problem_description}
\end{figure}

This is illustrated in Fig. \ref{fig:problem_description} for two initial points $y_1$ and $y_2$. Clearly, it requires deviating less from the deterministic flow when departing from $y_2$. However, it is necessary to analyse how \emph{likely} it is that the system starts from $y_2$ compared to $y_1$, according to $p_0$.

A solution to this problem is derived here using the notions described in the previous section. First, we put forward the following proposition which introduces a new variational problem. Later, we modify it to obtain an identical problem that can be numerically solved. The solution to the latter is the initial state and the trajectory of the system that exhibit highest probability density (in the sense of the large deviations principle) from which the system can depart and reach the boundary of the unsafe set $D$ at some time $t\in[0,T]$.
\begin{proposition}\label{proposition}
    Let $\phi_{y}^\ast(\cdot)$ be the maximum likelihood path hitting the boundary of $D$ at $t_y^\ast\leq T$, defining the quasi-potential $Q(y)=S_{t_y^\ast}(\phi_y^\ast)$, both as functions of the departing point $y$. \ra{Then, the most likely state of departure for small $\epsilon$ is the point $y^\ast$ minimizing the function
    \begin{equation}\label{eq:argmin_S_X0}
        \Gamma(y)=Q(y)+\epsilon S_0(y).\\[-3mm]
    \end{equation}}
\end{proposition}
\emph{Proof}. The proof of the previous proposition arises naturally from Bayes' theorem. \ra{For that, let the hitting event $E_y$ characterize a random variable $Z_y:\Omega\rightarrow\{0,1\}$, parameterized by $y$ and satisfying $Z_y(\omega)=1$ if $\omega\in E_y$, and $Z_y(\omega)=0$ otherwise. Since $E_y\in\mathcal{F}$ and $\{1\}$ is a Borel set, then $Z_y$ is measurable, and characterized by a probability mass function $p_{z|y}$ for $z\in\{0,1\}$ defined as
$$
p_{z|y}(z)=\bigl(\textbf{P}_y[\tau_D^\epsilon\leq T]\bigr)^z\bigl(1-\textbf{P}_y[\tau_D^\epsilon\leq T]\bigr)^{1-z}
$$
Consider as well the density $p_0$ of the initial (prior) distribution $\mu$ and its posterior distribution $\mu_\epsilon'$, updated by the observation of the hitting event, for which we assume the density $p_\epsilon(y|z)$. We now state Bayes' theorem for the densities of the random variables $Z_y$ and $X^\epsilon(0)$,
$$
    p_\epsilon(y|z)=\frac{p_{z|y}(z)p_0(y)}{p_{z}(z)}=\frac{p_{z|y}(z)p_0(y)}{\int p_{z|y}(z)p_0(y)dy}.
$$
Taking logarithms and multiplying by $\epsilon$, we get
$$
    \epsilon\ln p_\epsilon(y|z)=\epsilon\ln p_{z|y}(z)+\epsilon\ln p_0(y) - \epsilon\ln p_{z}(z).
$$
Let us focus only on the case $z=1$, and disregard the constant, non-zero term $p_z(1)=\textbf{P}_\mu[\tau_D^\epsilon\leq T]$,
$$
    \epsilon\ln p_\epsilon(y|z=1)\propto\epsilon\ln\textbf{P}_y[\tau_D^\epsilon\leq T] +\epsilon\ln p_0(y),
$$
where $\propto$ indicates a proportional relationship. Allow the process noise to be very small, i.e. $0<\epsilon\ll1$, enabling the notions of large deviations \eqref{eq:LDP} and \eqref{eq:quasipotential}.
$$
    \epsilon\ln p_\epsilon(y|z=1)\propto - Q(y) + o(1)-\epsilon S_0(y),
$$
where the small-$o$ notation denotes the small error associated to the LDP approximation.} It then holds that the argument $y^\ast$ that minimizes \ra{the right-hand side is the mode of the posterior distribution $p_\epsilon(y|z=1)$ in the low-noise approximation induced by the LDP, which yields
$$
   y^\ast = \underset{y}{\arg\max}~p_\epsilon(y|z=1)\overset{\epsilon}{\approx} \underset{y}{\arg\min}~\left(Q(y)\!+\!\epsilon S_0(y)\right),\\[-2mm]
$$
where the approximation stems from omitting the error terms in $o(1)$. Accuracy increases with decreasing $\epsilon$.}
\hfill$\square$

\begin{remark}
\rc{
Although the result is exact in the asymptotic limit $\epsilon\!\rightarrow\!0$, we instead consider small but finite $\epsilon$ that preserve the contribution of $S_0(y)$. From an engineering perspective, the LDP is therefore interpreted as a low-noise approximation rather than a strict equality.
}
\end{remark}

The previous proposition is instrumental to determine the worst initial condition $y^\ast$, in the probabilistic sense induced by the LDP. However, a closed-form solution of $\phi_{y}^\ast(\cdot)$ as a function of $y$ is unknown to us, hence preventing the direct use of \eqref{eq:argmin_S_X0}. Consider, then, the following.
\rc{
\begin{theorem}\label{thm:main}
For a system \eqref{eq:dSDE_X} with initial probability measure $\mu=p_0dy$, where $p_0$ is Gaussian as defined in \eqref{eq:mu}, the most probable path hitting $\partial D$ in the sense of Proposition \ref{proposition} is given by the pair $(\phi^\varoast,t^\varoast)$ where the minimum of the functional $J$ is attained:
\begin{equation}\label{eq:MAP}
\begin{aligned}
    Q' &\equiv \underset{t\leq T}{\textup{inf}}~\underset{\phi:\phi(t)\in\partial D}{\textup{min}}~J_t(\phi)\\ &= \underset{t\leq T}{\textup{inf}}~\underset{\phi:\phi(t)\in\partial D}{\textup{min}}\bigl\{S_t(\phi)+\epsilon S_0(\phi_0)\bigr\}.
    \end{aligned}
\end{equation}
\end{theorem}
}
\vspace{-5mm}
\emph{Proof}. We depart from Proposition \ref{proposition} and assume that the initial condition $y$ exhibiting highest posterior density is the one minimizing $\Gamma(y)$. Due to the lack of an analytic map $y\mapsto\phi_{y}^\ast(\cdot),t_y^\ast$, provided it requires to solve \eqref{eq:quasipotential}, then $y^\ast$, $\phi^\ast$ and  $t^\ast$ must be determined together. Let $\Gamma'$ be a function of these elements, defined as
$$\Gamma'(\phi,y,t)=S_t(\phi)+\epsilon S_0(y),\\[-3mm]$$
for which we consider the minimization
\begin{equation}
    \inf_{t\leq T}~ \inf_{y} ~\underset{\substack{\phi:\phi(0)=y,\\ \phi(t)\in\partial D}}{\textup{min}}~\Gamma'(\phi,y,t),
\end{equation}
and $(\phi^\ast,y^\ast,t^\ast)$ is the tuple where this minimum is attained. Note that the minimization with respect to $\phi$ is subject to the hard constraint $\phi(0)=y$, and later it becomes unconstrained in the minimization with respect to $y$, since we let $y\in\mathbb{R}^n$. Therefore, we can express the previous variational problem only on $\phi$ and $t$, with $\phi(0)$ unconstrained. This leads to the result in \eqref{eq:MAP}.
\hfill$\square$

As before, the involvement of $Q$ implies there is generally no unique solution. We also note that the function $S_0$ defines an additional cost on the initial state with a minimum at $x_0$. To distinguish the solution of this new variational problem from the solution of \eqref{eq:quasipotential}, we introduce the following definition.
\begin{definition}[Maximum a posteriori pathway]
    The extreme function $\phi^\varoast$, \rc{where the infimum $Q'$ is attained}, is the maximum a posteriori (MAP) pathway.
\end{definition}
\begin{remark}\label{rem:S0}
    The initial cost, $S_0$, in \eqref{eq:MAP} scales with the  noise $\epsilon$. As $\epsilon$ increases, the diffusion becomes more dispersed, increasing the likelihood of trajectories to reach $D$ while remaining closer to $x_0$. Thus, deviations from $x_0$ incur a higher effective cost, as captured in the VP.
\end{remark}
\begin{remark}
    Similarly to remark \ref{rem:S0}, the smaller the covariance $\Sigma$ in $S_0$, the less likely it is that sample paths depart further from $x_0$ will hit the unsafe set $D$
\end{remark}

\rc{Finally, we use \eqref{eq:LDP} to approximate the solution to the weak \emph{p}-safety problem \eqref{eq:psafety} in the low-noise case
\begin{equation}\label{eq:integral}
\textbf{P}_\mu[\tau_D^\epsilon\leq T]\overset{\epsilon}{\approx}\int\exp(-Q(y)/\epsilon)p_0(y)dy,\\[-1mm]
\end{equation}
with the caveat that it relies on a point-wise optimization problem for each $y$. We may use black-box or Bayesian optimization algorithms to solve \eqref{eq:integral} more efficiently; however, this lies outside the scope of this paper.}

\section{Application: Spacecraft collision}\label{sec:example}
The maximum a posteriori path is determined for a practical example\footnote{Code at \url{https://github.com/aitor-rg/LDT-safety}}: a collision between two space objects. The problem is solved numerically using CasADi, and the solution is cross-checked with the necessary conditions for optimality drawn from the maximum principle.

Consider two independent systems $(X_1(t))_{t\in I}$ and $(X_2(t))_{t\in I}$, each in $\mathbb{R}^6$, with unperturbed dynamics given by the gradient of the gravity potential $U$,\\[-0.5mm]
\begin{equation}
    \begin{aligned}
        dR_i(t) &= V_i(t) dt,\\
        dV_i(t) &= -\nabla U(R_i(t)) dt + \sqrt{\epsilon}\sigma dW_i(t),\\
    \end{aligned}
\end{equation}\\[0.5mm]
and initial conditions $X_i(0)\!\!\sim\!\!\{x_i,\Sigma_i\}$, for $i\in\{1,2\}$.

\begin{figure}[ht!]
    \centering
    \pgfplotstableread[col sep = comma]{figures/data/sat.dat}{\tabA}
\pgfplotstableread[col sep = comma]{figures/data/phi.dat}{\tabB}
\pgfplotstableread[col sep = comma]{figures/data/sat_end.dat}{\tabC}
\pgfplotstableread[col sep = comma]{figures/data/phi_end.dat}{\tabD}
\pgfplotstableread[col sep = comma]{figures/data/sat_start.dat}{\tabE}
\pgfplotstableread[col sep = comma]{figures/data/phi_start.dat}{\tabF}
\pgfplotstableread[col sep = comma]{figures/data/satphi_trans.dat}{\tabG}

\begin{tikzpicture}
\begin{groupplot}[
group style={
group name=my plots,
group size= 1 by 3,
horizontal sep=0cm,
vertical sep=-0.3cm,
},
footnotesize,
width=0.35\textwidth
]

\nextgroupplot[view={-25}{10}, axis lines=none,
                xmin = -3717800-500, xmax = -3717800+6500,
                ymin = 6114100-3500, ymax = 6114100+3500,
                zmin = 14944-3500, zmax = 14944+3500 ]
                
    \addplot3[black,very thick, smooth] table[x=s1x, y=s1y, z=s1z] {\tabE};
    \addplot3[red,very thick, dashed, smooth] table[x=s1x, y=s1y, z=s1z] {\tabF};
    \addplot3[black,very thick, smooth] table[x=s2x, y=s2y, z=s2z] {\tabG};
    \addplot3[red,very thick,dashed, smooth] table[x=p2x, y=p2y, z=p2z] {\tabG};
    \addplot3[black, only marks, mark=*, mark size=1.5] table {
        x   y   z
        -3.7179e+06 6.1141e+06  14944
        -3712900    6114100    14944
    } node[pos=0, below=1mm]{\footnotesize $r_1(0)$
    } node[pos=1, right=1mm]{\footnotesize $r_2(0)$};
    \addplot3[red, only marks, mark=*, mark size=1.5] table {
        x   y   z
        -3.7178e+06 6.1137e+06  16737        
        -3712700    6114200    13254
    } node[pos=0, right=1mm]{\footnotesize $\eta_1^\varoast(0)$
    } node[pos=1, below=1mm]{\footnotesize $\eta_2^\varoast(0)$};

\nextgroupplot[view={-20}{30}, grid = both, major grid style = {lightgray}, 
                minor grid style = {lightgray!25},minor tick num = 2,
                zmin = -6000000, zmax = 6000000, 
                xmin = -6000000, xmax = 6000000,
                ymin = -6000000, ymax = 6000000]

        \tikzset{mark options={draw=black,fill=black, mark size=2}}
        \addplot3[black,very thick, smooth] table[x=s1x, y=s1y, z=s1z] {\tabA}node[pos=0.75, left] {$r_1(t)$};
        \addplot3[black,very thick, smooth] table[x=s2x, y=s2y, z=s2z] {\tabA} node[pos=0.75, right] {$r_2(t)$};
        \addplot3[red, very thick, smooth, dashed] table[x=s1x, y=s1y, z=s1z] {\tabB};
        \addplot3[red, very thick, smooth, dashed] table[x=s2x, y=s2y, z=s2z] {\tabB};
        \addplot3[black, only marks, mark=o, mark size=2] table {
            x   y   z
            -3.7179e+06 6.1141e+06  14944
            6.4654e+06  3.9334e+06  -4.9558e+05        
        };

        \addplot3[opacity=0.3,black, only marks, mark=*, mark options={mark size=18}] table {
            x   y   z
            0   0   0
        };
        
\nextgroupplot[view={-20}{20}, axis lines=none]
        \addplot3[black,very thick, smooth] table[x=s1x, y=s1y, z=s1z] {\tabC} node[pos=0.2, below=1mm]{$r_1(t)$};
        \addplot3[black,very thick, smooth] table[x=s2x, y=s2y, z=s2z] {\tabC} node[pos=0.8, above=1mm]{$r_2(t)$};
        \addplot3[red,  very thick, smooth, dashed] table[x=s1x, y=s1y, z=s1z] {\tabD} node[pos=0.5, red, above=1mm]{$\eta_1^\varoast(t)$};
        \addplot3[red,  very thick, smooth, dashed] table[x=s2x, y=s2y, z=s2z] {\tabD} node[pos=0.5, red, below=0.5mm]{$\eta_2^\varoast(t)$};
        \addplot3[black, only marks, mark=o, mark options={mark size=2}] table {
            x   y   z
            0.3739e+06  0.2501e+06  -7.1410e+06 
        } node[pos=0, above right=0mm]{$D$};
\end{groupplot}

\newcommand\s{0.1}
\draw[black] (0.76,-6.18) rectangle ++(3.87,1.92); 
\draw[black] (2.23,-2.97) rectangle ++(3.85*\s,1.98*\s); 
\draw[black] (2.39,-2.97) -- (2.65,-4.24);

\draw[black] (0.3,1) rectangle ++(1.9,2.6); 
\draw[black] (0.53,-1.24) rectangle ++(1.9*\s, 2.6*\s); 
\draw[black] (0.53,-0.98) -- (0.3,1);

\draw[black] (2.5,0.5) rectangle ++(1.9,3.22); 
\draw[black] (3.64,-1.15) rectangle ++(1.9*\s, 2.6*\s); 
\draw[black] (3.64,-1.15) -- (2.5,0.5);
    
\end{tikzpicture}
    \vspace{-7mm}
    \caption{Conjunction geometry and MAP path solution.}
    \label{fig:geometry}
\end{figure}
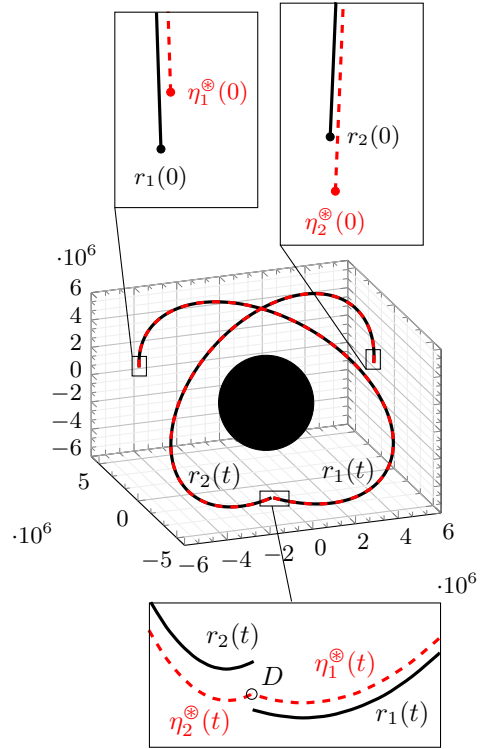

For convenience, we assume that both objects have the same noise parameter $\epsilon$ and constant diffusion $\sigma$ with $a_i=\sigma\sigma^\top$, and we augment the state vector $X(t)\equiv(X_1(t),X_2(t))$ so that $x_0\equiv(x_1,x_2)$ with $x_i=(r_i(t),v_i(t))$, $\Sigma\equiv\text{diag}(\Sigma_1,\Sigma_2)$ and $\phi(t)\equiv(\phi_1(t),\phi_2(t))$ with $\phi_i(t)=(\eta_i(t),\nu_i(t))$ for $i\in\{1,2\}$, as presented in section \ref{sec:preliminaries}. Additionally, let $\eta(t)\equiv(\eta_1(t),\eta_2(t))$ and $\nu(t)\equiv(\nu_1(t),\nu_2(t))$. The mean initial conditions are defined in such a way that both space objects come to a distance of 7km from each other, as illustrated in Fig. \ref{fig:geometry} where the deterministic trajectory (continuous black) and the MAP path solution (dashed red) are shown.

The complete variational problem is\footnote{Notice that solving for the path $\phi$ and solving for the pair $(\phi_0,w)$ as an optimal control problem is equivalent given the differential equation \eqref{eq:phi_dynamics_mech}.}:
\begin{equation}
    \begin{aligned}
        \underset{\phi_0,w,T}{\text{min}}&~\left\{\int_0^T\tfrac{1}{2}||w||^2 dt+ \tfrac{\epsilon}{2}||\phi_0-x_0||_\Sigma^2\right\}\\[1mm]
        \text{s.t.}\hspace{0mm} & \hspace{5mm}\dot{\eta}=\nu,\\[1mm] & \hspace{5mm}\dot{\nu}=-\nabla U(\eta) +\sigma w, \quad\forall t\in [0,T],\\[1mm]
        & \hspace{5mm}\phi(0) \text{ free}, \ \eta(T)\in \partial \tilde{D},~ T\in[T_1,T_2]\\
    \end{aligned}
\end{equation}
where $\nabla U(\eta)\equiv(\nabla U(\eta_1),\nabla U(\eta_2))$, $a \equiv\text{diag}(a_1,a_2)$ and $ \tilde{D}\equiv\{z\in\mathbb{R}^{12}:f(z)\leq0\}$. The function $f$ is a terminal constraint defining a collision configuration, e.g,
$$f(\phi(T)) = ||\eta_1(T)-\eta_2(T)||^2- \gamma \leq0,$$
 with the parameter $\gamma\geq0$ being a safety critical relative distance between the two objects.

 The \emph{Hamiltonian} of the variational problem can be defined for an adjoint state $\lambda(t):=(\lambda_\eta(t),\lambda_\nu(t))$ as
\begin{equation*}
    H(\phi,w,\lambda)=-\frac{1}{2}w^\top w + \begin{bmatrix}
        \lambda_\eta\\
        \lambda_\nu
    \end{bmatrix}^\top\begin{bmatrix}
        \nu\\
        -\nabla U(\eta)+\sigma w
    \end{bmatrix},\\[-1mm]
\end{equation*}
which provides the MAP path and adjoint equations
\begin{subequations}\label{eq:MAPeq}
    \begin{alignat}{4}
        \dot{\phi}^\varoast & =   \nabla_{\!\lambda} H = \begin{bmatrix}
            \nu^\varoast\\
            -\nabla U(\eta^\varoast) + a \lambda_\nu
        \end{bmatrix},\\
        \dot{\lambda}& = -\nabla_{\!\phi} H = \begin{bmatrix}
            0 & \nabla^2 U(\eta^\varoast)\\
            -I & 0
        \end{bmatrix}\lambda,\label{eq:map_nu}
    \end{alignat}\\[-4mm]
\end{subequations}
since the optimal deviation is $w^\varoast=\sigma^\top\lambda_\nu$ given by the \emph{maximizing condition} of the maximum principle. This is the optimal deviation from the deterministic vector field satisfying $f(\phi(T))\leq0$. If the trajectory of the deterministic system already satisfy this inequality constraint, then the optimal solution would yield $w=\lambda=0$, $\forall t\in I$.

The transversality conditions \cite{Atle:87} are, with $\alpha\geq0$
\begin{subequations}\label{eq:trans_cond}
    \begin{alignat}{4}
        \lambda(0) & = -\epsilon\nabla_{\!\phi} S_0(\phi^\varoast(0)), \label{eq:new_trans_cond}\\ 
        \lambda(T^\varoast)  & =  \alpha\nabla_{\!\phi} f (\phi^\varoast(T^\varoast)),\label{eq:final_trans_cond}\\ 
        \alpha f(\phi(T^\varoast)) & = 0,
    \end{alignat}\\[-5mm]
\end{subequations}
for the adjoint state, and the temporal condition
\begin{equation}\label{eq:condT}
    H(\phi^\varoast(T^\varoast),w^\varoast(T^\varoast),\lambda(T^\varoast)) = 0.
\end{equation}
Note that the ML and MAP path problems differ slightly in the necessary conditions for optimality. The hard constraint $\phi_0=x_0$ in the ML path problem is replaced by a transversality condition on the co-state $\lambda$ of the MAP path problem, that is equation \eqref{eq:new_trans_cond}, leaving $\phi_0$ free.

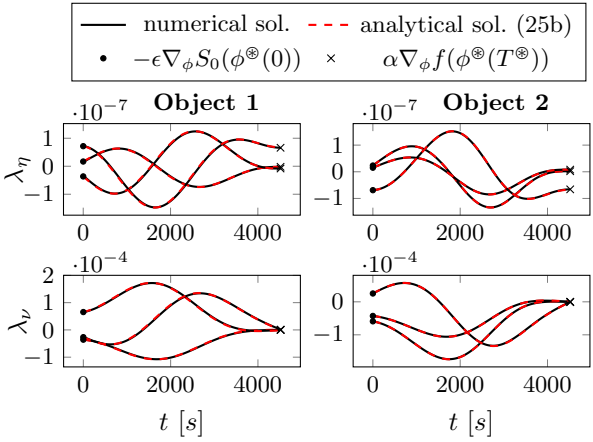
\begin{figure}[ht!]
    \begin{center}~~~~\ref{leg}\end{center}
    \begin{flushleft}
    \pgfplotstableread[col sep = comma]{figures/data/lam1_num.dat}{\condA}
\pgfplotstableread[col sep = comma]{figures/data/lam2_num.dat}{\condB}
\pgfplotstableread[col sep = comma]{figures/data/lam1_dyn.dat}{\condC}
\pgfplotstableread[col sep = comma]{figures/data/lam2_dyn.dat}{\condD}

\begin{tikzpicture}
\begin{groupplot}[
        group style={
            group name=my plots,
            group size= 2 by 2,
            horizontal sep=0.7cm,
            vertical sep=0.8cm,
        },
        footnotesize,
        width=0.267\textwidth,
        height=2.8cm,
        ylabel near ticks,
        xlabel near ticks,
        xtick={0,100,200},
        xticklabels={$0$,$2000$,$4000$},
        legend columns=2,
        legend entries={numerical sol.~, analytical sol. \eqref{eq:map_nu}, $-\epsilon\nabla_\phi S_0(\phi^\varoast(0))$, $\alpha\nabla_\phi f(\phi^\varoast(T^\varoast))$},
        legend to name=leg,
        legend style={draw opacity = 1}
        ]
        
\nextgroupplot[title={~~~~~~\textbf{Object 1}},ylabel near ticks,
    xlabel near ticks, ylabel={$\lambda_\eta$}, y label style={at={(axis description cs:-0.1,0.5)}}]
    \addplot[black,thick, smooth] table[y=l1x, x=t] {\condA};
    \addplot[red,dashed,thick, smooth] table[y=l1x, x=t] {\condC};
    \addplot[black,only marks, mark=*, mark size=1] table {
            x   y
            0   1.6843e-08
            0   -3.6397e-08 
            0   7.1668e-08
        };
    \addplot[black,only marks, mark=x, mark size=2] table {
            x   y
            226   -2.0453e-09
            226   -8.1743e-09
            226   6.5939e-08
        };
    \addplot[black,thick, smooth] table[y=l1y, x=t] {\condA};
    \addplot[red,dashed,thick, smooth] table[y=l1y, x=t] {\condC};
    \addplot[black,thick, smooth] table[y=l1z, x=t] {\condA};
    \addplot[red,dashed,thick, smooth] table[y=l1z, x=t] {\condC};
    
\nextgroupplot[title={~~~~~~\textbf{Object 2}}]
    \addplot[black,thick, smooth] table[y=l1x, x=t] {\condB};
    \addplot[red,dashed,thick, smooth] table[y=l1x, x=t] {\condD};
    \addplot[black,only marks, mark=*, mark size=1] table {
            x   y
            0   2.3352e-08
            0   1.5245e-08
            0   -6.9002e-08
        };
    \addplot[black,only marks, mark=x, mark size=2] table {
            x   y
            226   2.0448e-09
            226   8.1744e-09
            226   -6.5939e-08
        };
    \addplot[black,thick, smooth] table[y=l1y, x=t] {\condB};
    \addplot[red,dashed,thick, smooth] table[y=l1y, x=t] {\condD};
    \addplot[black,thick, smooth] table[y=l1z, x=t] {\condB};
    \addplot[red,dashed,thick, smooth] table[y=l1z, x=t] {\condD};
    
\nextgroupplot[ylabel near ticks,
    xlabel near ticks,xlabel={$t~[s]$}, ylabel={$\lambda_\nu$}, y label style={at={(axis description cs:-0.1,0.5)}}]
    \addplot[black,thick, smooth] table[y=l1vx, x=t] {\condA};
    \addplot[red,dashed,thick, smooth] table[y=l1vx, x=t] {\condC};
    \addplot[black,only marks, mark=*, mark size=1] table {
            x   y
            0   -3.5785e-05
            0   6.5496e-05
            0   -2.6657e-05
        };
    \addplot[black,only marks, mark=x, mark size=2] table {
            x   y
            226   -3.4751e-08
            226   -1.3897e-07
            226   1.1207e-06
        };
    \addplot[black,thick, smooth] table[y=l1vy, x=t] {\condA};
    \addplot[red,dashed,thick, smooth] table[y=l1vy, x=t] {\condC};
    \addplot[black,thick, smooth] table[y=l1vz, x=t] {\condA};
    \addplot[red,dashed,thick, smooth] table[y=l1vz, x=t] {\condC};
    
\nextgroupplot[ylabel near ticks,
    xlabel near ticks,xlabel={$t~[s]$}]
    \addplot[black,thick, smooth] table[y=l1vx, x=t] {\condB};
    \addplot[red,dashed,thick, smooth] table[y=l1vx, x=t] {\condD};
    \addplot[black,only marks, mark=*, mark size=1] table {
            x   y
            0   -5.9107e-05
            0   -4.282e-05
            0   2.5271e-05
        };
    \addplot[black,only marks, mark=x, mark size=2] table {
            x   y
            226   3.4749e-08
            226   1.3897e-07
            226   -1.1207e-06
        };
    \addplot[black,thick, smooth] table[y=l1vy, x=t] {\condB};
    \addplot[red,dashed,thick, smooth] table[y=l1vy, x=t] {\condD};
    \addplot[black,thick, smooth] table[y=l1vz, x=t] {\condB};
    \addplot[red,dashed,thick, smooth] table[y=l1vz, x=t] {\condD};
    
\end{groupplot}
\end{tikzpicture}
    \end{flushleft}
    \vspace{-2mm}
    \caption{Numerical solution of $\lambda$ superposed with the solution of (29b), and boundary conditions (30a) and (30b).}
    \label{fig:transversality}
\end{figure}

Conditions \eqref{eq:trans_cond} and \eqref{eq:condT} are satisfied by the obtained numerical solution, which implies that $\phi^\varoast$ is an extremal candidate. The initial co-state $\lambda(0)$ obtained from the numerical solver has been used to integrate the differential equation \eqref{eq:map_nu} and compare it to the trajectory of the entire numerical solution, as shown in Fig. \ref{fig:transversality}. The boundary conditions \eqref{eq:trans_cond} are also shown. The optimal final time is found to be $T^\varoast=4518.1$ seconds and $\alpha=6.81\cdot10^{-10}$.
\vspace{-2mm}

\begin{figure}[ht!]
    \centering
    \begin{center}\ref{leg2}~~\end{center}
    \begin{flushleft}
    \pgfplotstableread[col sep = comma]{figures/data/action_dagger.dat}{\devA}
\pgfplotstableread[col sep = comma]{figures/data/action_star.dat}{\devB}
\pgfplotstableread[col sep = comma]{figures/data/w_dagger.dat}{\devC}
\pgfplotstableread[col sep = comma]{figures/data/w_star.dat}{\devD}

\begin{tikzpicture}
\begin{groupplot}[
            group style={
            group name=my plots,
            group size= 2 by 1,
            horizontal sep=1.1cm,
            vertical sep=1cm,
        },
        footnotesize,
        width=0.25\textwidth,
        height=2.8cm,
        xtick={0,100,200},
        xticklabels={$0$,$2000$,$4000$},
        legend columns=2,
        legend entries={ML,MAP},
        legend to name=leg2,
        ]

\nextgroupplot[xlabel={$t~[s]$},ylabel={$S_T$},y label style={at={(axis description cs:0.2,0.5)}}]
    \addplot[black,very thick, smooth] table[x=t, y=y1]{\devB};
    \addplot[black,very thick, dashed, smooth] table[x=t, y=y1]{\devA};
    
\nextgroupplot[xlabel={$t~[s]$},ylabel={$||w||$},y label style={at={(axis description cs:0.2,0.5)}}]
    \addplot[black,very thick, smooth] table[x=t, y=y1]{\devD};
    \addplot[black,very thick, dashed, smooth] table[x=t, y=y1]{\devC};

\end{groupplot}
\end{tikzpicture}
    \vspace{-2mm}
    \end{flushleft}
    \caption{Comparison between the action $S_T$ and magnitude of the deviation $||w||$ exerted by the ML and MAP paths.}
    \label{fig:deviations}
\end{figure}
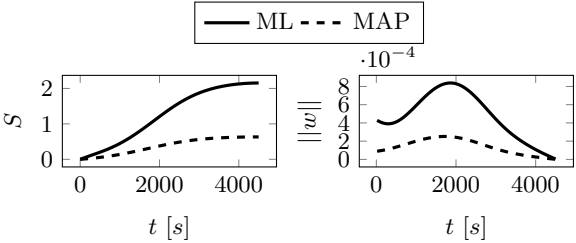

Finally, in Fig. \!\ref{fig:deviations} we compare the action $S_T$ and the magnitude of the deviation $w$ at each time instant $t\in[0,T^\varoast]$ for both, the ML and MAP path solutions. As expected, the MAP path requires less action, and thus smaller deviations, than the ML path to realize the URE.
\vspace{-3mm}

\section{Conclusions}\label{sec:conclusions}
We have introduced a new tool to analyze safety-critical systems with uncertain initial conditions in scenarios where unsafe events are rare, which challenge the numerical computation of hitting probabilities via standard methods. \rc{Combining Bayes' theorem and a principle of Large Deviations, we have identified the posterior density of the initial states in the low-noise case; this expression can be maximized to approximate the most probable initial conditions, paths and hitting times in which the rare events are observed}. The solutions enable complementary analyses, e.g. the study of posterior events (in our example, the debris dispersion generated after the collision, which is highly dependent on the collision geometry) or the design of maneuvers that avoid these paths. \rc{We also find a computational expression for the weak \emph{p}-safety problem, which can be effectively solved through black-box or Bayesian optimization algorithms.} Finally, the feasibility of the proposed procedure has been tested under a real high-dimensional and non-linear use-case. The results are obtained numerically, and generally has no unique solutions. Heuristically, however, one could solve the proposed variational problem multiple times with different initial guesses in order to collect different local minimizers. The collected extreme paths can be used to design importance sampling algorithms or control laws that steer the deterministic trajectory of the system away from these minimizers.

\bibliographystyle{plain}        
\bibliography{autosam}           



\end{document}